\documentclass[12pt]{amsart}
\usepackage{amsmath,amssymb,amsfonts,amsthm,amstext,verbatim,url}
\usepackage{graphicx}
\RequirePackage[colorlinks,citecolor=blue,urlcolor=blue]{hyperref}

\theoremstyle{plain}

\newcounter{qn}

\newtheorem{theorem}{Theorem}
\newtheorem{definition}[theorem]{Definition}

\newtheorem{example}[theorem]{Example}

\newtheorem{question}[qn]{Question}
\numberwithin{equation}{section}
\numberwithin{theorem}{section}
\numberwithin{figure}{section}

\newcounter{mycount}
\newenvironment{romlist}{\begin{list}{\rm(\roman{mycount})}%
   {\usecounter{mycount}\labelwidth=1cm\itemsep 0pt}}{\end{list}}
\newenvironment{numlist}{\begin{list}{\rm{\arabic{mycount}}.}%
   {\usecounter{mycount}\labelwidth=1cm\itemsep 0pt}}{\end{list}}
\newenvironment{letlist}{\begin{list}{\rm(\alph{mycount})}%
   {\usecounter{mycount}\labelwidth=1cm\itemsep 0pt}}{\end{list}}

\newcommand\s{\sigma}

\newcommand\ol{\overline}
\newcommand\oo{\infty}
\newcommand\De{\Delta}
\newcommand\LL{{\mathbb L}}
\newcommand\HH{{\mathbb H}}
\newcommand\NN{{\mathbb N}}
\newcommand\TT{{\mathbb T}}

\newcommand\sP{{\mathcal P}}

\newcommand\sA{{\mathcal A}}

\newcommand\sN{{\mathcal N}}

\newcommand\ZZ{{\mathbb Z}}
\newcommand\RR{{\mathbb R}}
\newcommand\BB{{\mathbb B}}

\newcommand\wt{\widetilde}
\newcommand\om{\omega}

\renewcommand\a{\alpha}
\newcommand\Ga{\Gamma}

\newcommand\eps{\epsilon}

\newcommand\g{\gamma}

\newcommand\fish{F}

\newcommand\olG{{\ol G}}
\newcommand\olV{{\ol V}}
\newcommand\olE{{\ol E}}

\newcommand\vG{{\vec G}}

\newcommand\vE{{\vec E}}

\newcommand\vpi{\vec\pi}
\newcommand\vmu{\vec\mu}

\newcommand\Aut{\text{\rm Aut}}
\newcommand\pd{\partial}

\newcommand\id{\iota}

\newcommand\fcp{finite coset property}
\newcommand\SLE{\text{\rm SLE}}

\newcommand\muF{\mu_{\text{\rm F}}}
\newcommand\siF{\s^{\text{\rm F}}}
\newcommand\pc{p_{\text{\rm c}}}
\newcommand\Tc{T_{\text{\rm c}}}
\newcommand\ollG{G'}
\newcommand\ollE{E'}

\begin{document}
\title[Counting self-avoiding walks]
{Counting self-avoiding walks}
\author[Grimmett]{Geoffrey R.\ Grimmett}
\address{Statistical Laboratory, Centre for
Mathematical Sciences, Cambridge University, Wilberforce Road,
Cambridge CB3 0WB, UK} 
\email{\{g.r.grimmett, z.li\}@statslab.cam.ac.uk}
\urladdr{\url{http://www.statslab.cam.ac.uk/~grg/}}
\urladdr{\url{http://www.statslab.cam.ac.uk/~zl296/}}

\author[Li]{Zhongyang Li}
\begin{abstract}
The \emph{connective constant} $\mu(G)$
of a graph $G$ is the asymptotic growth
rate of the number of self-avoiding walks on $G$
from a given starting vertex.
We survey three aspects of the dependence of the connective constant
on the underlying graph $G$. Firstly, when $G$ is cubic,
we study the effect on $\mu(G)$ of the Fisher transformation 
(that is, the replacement of vertices by triangles).
Secondly, we discuss upper and lower bounds for $\mu(G)$
when $G$ is regular.  Thirdly, we present strict inequalities
for the connective constants $\mu(G)$ of vertex-transitive graphs $G$, as $G$ varies. 
As a consequence of the last, the connective constant of a Cayley graph 
of a finitely generated group decreases strictly when
a new relator is added, and increases strictly when a 
non-trivial group element is declared to be a generator. 
Special prominence is given to open problems.
\end{abstract}

\date{19 April 2013, revised 20 March 2015}  

\maketitle

\section{introduction}
\begin{figure}[htbp]
\centerline{\includegraphics[width=0.3\textwidth]{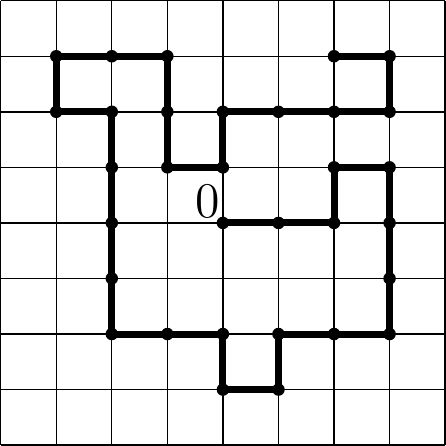}}
   \caption{A SAW from the origin of the square lattice.}
\label{fig:saws}
\end{figure}

A \emph{self-avoiding walk} (abbreviated to SAW) on a 
graph $G=(V,E)$ is a path that visits no vertex
more than once. An example of a  SAW on the square 
lattice is drawn in Figure \ref{fig:saws}.
SAWs were first introduced in the chemical theory of polymerization (see Flory \cite{f}),
and their critical behaviour has attracted the 
abundant attention since of mathematicians and physicists 
(see, for example, the book of Madras and Slade \cite{ms}
or the lecture notes \cite{bdgs}). 

Let $\s_n(v)$
be the number of $n$-step SAWs on $G$ starting  at the vertex $v$.
The following fundamental theorem of Hammersley asserts the existence of an
asymptotic growth rate for $\s_n(v)$ as $n \to \oo$. 
(See the start of Section \ref{sec:fisher} for a 
definition of (quasi-)transitivity.)

\begin{theorem}\cite{jmhII}\label{jmh}
Let $G=(V,E)$ be an infinite, connected, quasi-transitive graph with finite vertex-degrees. There exists
$\mu=\mu(G)\in[1,\oo)$, called the \emph{connective constant} of $G$,  such that
\begin{equation}\label{connconst}
\lim_{n \to \oo} \s_n(v)^{1/n} = \mu, \qquad v \in V.
\end{equation}
\end{theorem}

We review here recent work from \cite{GrL1,GrL2,GrL3} on the dependence
of $\mu(G)$ on the choice of graph $G$. 

\begin{figure}[htb]
 \centering
\includegraphics[width=0.8\textwidth]{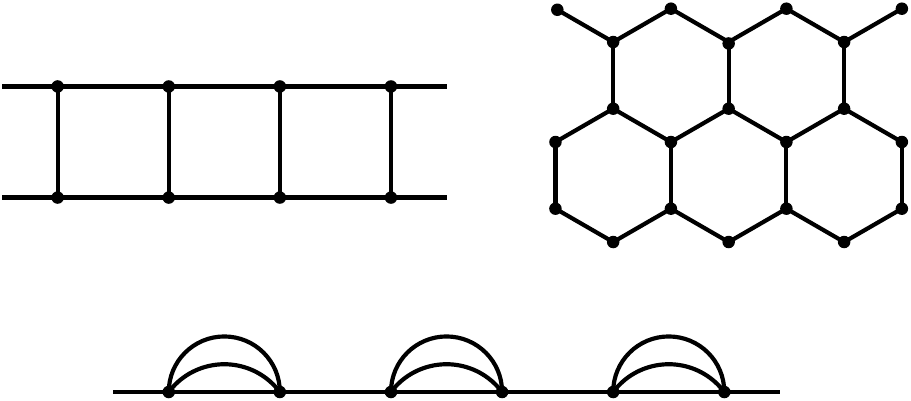}
  \caption{Three regular graphs: the ladder graph $\LL$,  
the hexagonal tiling $\HH$
of the plane, and the bridge graph $\BB_\Delta$ 
(with $\Delta=4$) obtained from $\ZZ$ by joining
every alternate pair of consecutive vertices by $\Delta-1$ parallel edges.}
  \label{fig:ladder-hex}
\end{figure}

For what graphs $G$ is $\mu(G)$ known exactly? 
There are a number of such graphs, 
which should be regarded as atypical in this regard.
We mention the \emph{ladder} $\LL$, the hexagonal lattice $\HH$, and the 
\emph{bridge graph}  $\BB_\De$ with degree $\De\ge 2$ of Figure \ref{fig:ladder-hex}, for which
\begin{equation}\label{2}
\mu(\LL) = \tfrac12(1 + \sqrt 5), \quad \mu(\HH) = \sqrt{2+\sqrt 2},
\quad \mu(\BB_\De) = \sqrt {\De-1}.
\end{equation}
See \cite[p.\ 184]{AJ90} and \cite{ds}
for the first two calculations. 

In contrast, the
connective constant of the square grid $\ZZ^2$ is unknown, and a substantial
amount of work has been devoted to obtaining good bounds. 
The best rigorous bounds known currently to the authors are those of \cite{j04,PT}, namely
(to $5$ significant figures)
\begin{equation*}
2.6256\leq \mu(\ZZ^2)\leq 2.6792,
\end{equation*}
and more precise numerical estimates are available, including the estimate
 $\mu \approx 2.63815\dots$ of \cite{Jens03}.

We make some remarks about the three graphs of Figure \ref{fig:ladder-hex}.
There is a correspondence between the Fibonacci sequence 
and counts of SAWs on the ladder graph $\LL$, whereby one obtains that 
$\mu(\LL)$ equals the golden ratio $\phi := \frac12(1+\sqrt 5)$.
We ask in Section \ref{sec:lb} (of the current article)
whether $\mu(G) \ge \phi$ for all simple, cubic, vertex-transitive graphs.
Amongst a certain category of $\De$-regular 
graphs permitted to possess multiple edges, the bridge graph
$\BB_\De$ is extremal in the sense that $\mu(\BB_\De)$ is smallest.
See the discussion of Section \ref{sec:lb}.

The proof that $\mu(\HH)=\sqrt{2+\sqrt{2}}$ by 
Duminil-Copin and Smirnov \cite{ds}
is a very significant recent result. The value $\sqrt{2+\sqrt 2}$ 
emerged in the physics literature through work  of Nienhuis \cite{Nien}
motivated originally by 
renormalization group theory. Its proof in \cite{ds} is based
on the construction of an observable with some properties of discrete holomorphicity,
complemented by a neat use of the bridge decomposition introduced by Hammersley and Welsh
\cite{HW62}.

It is a beautiful open problem to prove that
a random $n$-step SAW from the origin of $\ZZ^2$ converges, when suitably re-scaled,
to the Schramm--Loewner curve $\SLE_{8/3}$.
This important conjecture has been discussed and formalized 
by Lawler, Schramm, and Werner \cite{lsw}. 

\begin{question}
Does a uniformly distributed  $n$-step SAW on $\ZZ^2$ converge,
when suitably rescaled,
to the random curve $\SLE_{8/3}$?
\end{question}

There is an important class of results usually referred to as the `pattern theorem'.
In Kesten's original paper \cite{hkI} devoted to $\ZZ^2$,
a \emph{proper internal pattern} $\sP$  is defined as a finite SAW with the property
that, for any $k \ge 1$, there exists a SAW containing at least $k$ translates
of $\sP$. The pattern theorem states that: for a given proper internal pattern
$\sP$, there exists $a>0$ such that the number of $n$-step SAWs 
from the origin $0$,  
containing no more than $an$ translates of $\sP$, 
is exponentially smaller than the total $\s_n := \s_n(0)$. 

The pattern theorem may be used to prove for this bipartite graph
that
\begin{equation*}
\lim_{n\rightarrow\infty}\frac{\sigma_{n+2}}{\sigma_n}=\mu^2.
\end{equation*} 
The following stronger statement has been open since Kesten's paper \cite{hkI},
see the discussion at \cite[p.\ 244]{ms}.

\begin{question}
Is it the case for SAWs on $\ZZ^2$ that $\s_{n+1}/\s_n \to \mu$?
\end{question}

Previous work on SAWs tends to have been focussed on specific graphs such as 
the cubic lattices $\ZZ^d$ and certain two-dimensional lattices. 
In contrast, the results of \cite{GrL1,GrL2,GrL3}, surveyed here, 
are directed at
large classes of regular graphs, often but not 
exclusively \emph{transitive} graphs.
The work reviewed here may be the first systematic study of 
SAWs on general vertex-transitive and 
quasi-transitive graphs.

In Section \ref{sec:fisher}, we describe the Fisher transformation
and its effect on counts of SAWs for cubic (or semi-cubic) graphs. 
Universal bounds for the connective constant 
$\mu(G)$ of $\De$-regular graphs $G$ are presented in
Section \ref{sec:lb}. 
Strict inequalities for $\mu(G)$ are presented in Section \ref{sec:si}.
The question addressed there is the following: if $G$ is a strict subgraph of $G'$, under
what conditions on the pair $(G,G')$ is it the case
that $\mu(G) < \mu(G')$? Two sufficient conditions of an algebraic
nature are presented, with applications (in Section \ref{sec:cayley})
to Cayley graphs of
finitely generated groups.

A number of questions are included in this review. 
The inclusion of a question does not of itself imply
either difficulty or importance. 

\smallskip
\noindent
\emph{Note added at revision.} Since this review was written in 2013, the authors
have continued their project with papers \cite{GrL4,GrL5},
which are directed at the question of `locality' of connective constants: to what degree
is the value of the connective constant of a vertex-transitive graph determined by knowledge of
a large ball centred at a given vertex?

\section{The Fisher transformation and the golden mean}\label{sec:fisher}

This section is devoted to a summary of the effect on $\mu(G)$ of
the so-called \emph{Fisher transformation}. We begin with a discussion of transitivity.

The automorphism group of the graph $G=(V,E)$ is denoted $\Aut(G)$,
and the identity automorphism is written $\id$.
A subgroup $\Ga\subseteq\Aut(G)$ is said to \emph{act transitively} on 
$G$ if, for $v,w\in V$,
there exists $\g\in\Ga$ with $\g v=w$. It is said to \emph{act quasi-transitively} 
if there exists a finite set $W$
of vertices such that, for $v\in V$, there exist 
$w\in W$ and $\g\in \Ga$ with
$\g v=w$. The graph is called \emph{vertex-transitive} 
(respectively, \emph{quasi-transitive}) 
if $\Aut(G)$ acts transitively
(respectively, quasi-transitively) on $G$.

An automorphism $\g$ is said to \emph{fix} a vertex $v$ if $\g v = v$.
The subgroup $\Ga$ is said to \emph{act freely} on $G$ (or on the vertex-set $V$) if:
whenever there exist $\g \in \Ga$ and $v\in V$ with $\g v = v$, then $\g = \id$.

\begin{figure}[htb]
 \centering
    \includegraphics[width=0.5\textwidth]{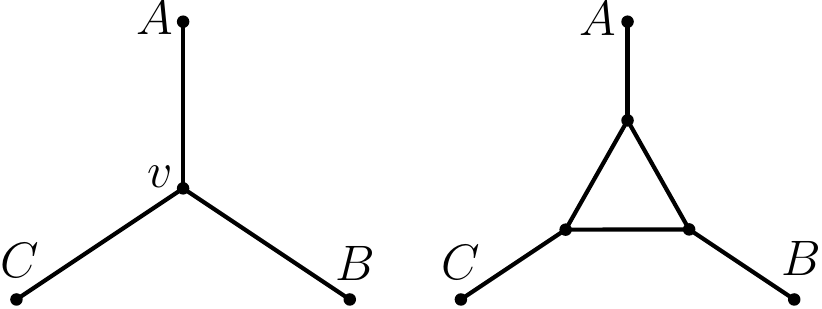}
  \caption{The Fisher triangulation of the star.}
  \label{fig:fisher}
\end{figure}

Let $v$ be a vertex with degree $3$. The \emph{Fisher transformation} acts at $v$ by 
replacing it by a triangle, as illustrated in 
Figure \ref{fig:fisher}. 
The Fisher transformation has been valuable in the study of the relations 
between Ising, dimer, and general vertex models
(see \cite{bdet,fisher,zli,li}), and also
in the calculation of the connective constant of the Archimedean lattice $(3,12^2)$ (see,
for example, \cite{g,GPR,jg}).
The Fisher transformation may be applied at every vertex of a cubic 
graph (that is, a graph with every vertex of degree $3$),
of which the hexagonal and square/octagon lattices are examples.
We describe next the Fisher transformation in the 
context of self-avoiding walks.

A graph $G$ is called \emph{simple} if it has no multiple edges.
Assume
that $G=(V,E)$ is quasi-transitive, connected, and simple. By Theorem \ref{jmh}, $G$ has
a well-defined connective constant $\mu=\mu(G)$ satisfying \eqref{connconst}.
Suppose, in addition, that $G$ is cubic, and
write $F(G)$ for the graph obtained by applying the Fisher transformation
at every vertex. The automorphism group of $G$ induces an automorphism subgroup of $F(G)$,
so that $F(G)$ is quasi-transitive and has a well-defined connective constant.
It is noted in \cite{GrL2}, and probably elsewhere also, that the connective constants of
$G$ and $F(G)$ have a simple relationship. This conclusion, and its iteration, 
are given in the
next theorem, in which $\phi := \frac12(1+\sqrt 5)$ denotes the golden mean.

\begin{theorem} \cite[Thm 3.1]{GrL2} \label{thm:main2}
Let $G$ be an infinite, quasi-transitive, connected, cubic graph, and consider the
sequence $(G_k: k=0,1,2,\dots)$ given by $G_0=G$ and $G_{k+1} = \fish(G_k)$.
\begin{letlist}
\item The connective constants $\mu_k := \mu(G_k)$ satisfy
$\mu_k^{-1} = g(\mu_{k+1}^{-1})$ where
$g(x)= x^2 + x^3$.
\item 
The sequence $\mu_k$ converges monotonely to $\phi$,
and 
\begin{equation*}
- \left(\tfrac 47\right)^k \le \mu_k^{-1} - \phi^{-1} 
\le \bigl[\tfrac12(7-\sqrt 5)\bigr]^{-k}, \qquad k \ge 1.
\end{equation*}
\end{letlist}
\end{theorem}

The idea underlying part (a) is that, at each vertex $v$ visited by a SAW $\pi$ on $G_k$, one 
may replace that vertex by either of the two paths around the `Fisher triangle'
at $v$. Some book-keeping is necessary with this argument, and this is best done 
via generating functions \eqref{eq:gf}. 

We turn now to the Fisher transformation in the context of a `semi-cubic' graph.
A graph is called \emph{bipartite} if its vertices can be coloured black or white 
in such a way that every edge links a black vertex and a white vertex.

\begin{theorem} \cite[Thm 3.3]{GrL2} \label{sf}
Let $G$ be an infinite, connected, bipartite graph with vertex-sets coloured black
and white, and suppose the coloured graph is quasi-transitive, and every black vertex has degree $3$. Let
$\wt{G}$ be the graph obtained by applying the Fisher transformation at each black vertex. 
The connective 
constants $\mu$ and $\wt\mu$ of $G$ and $\wt{G}$, respectively, 
satisfy $\mu^{-2}=h(\wt{\mu}^{-1})$, where $h(x)=x^3+x^4$.
\end{theorem}

\begin{figure}[htb]
 \centering
    \includegraphics[width=0.5\textwidth]{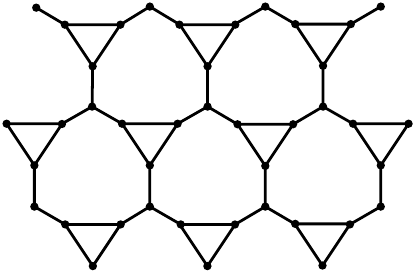}
  \caption{The lattice $\wt\HH$ derived from the hexagonal lattice $\HH$ by
applying the Fisher transformation at alternate vertices. Its connective constant
$\wt\mu$ is the root of the equation
$x^{-3}+x^{-4} = 1/(2+\sqrt 2)$.}
  \label{fig:newlatt}
\end{figure}

\begin{example}
Theorem \ref{sf} implies an exact value of a connective constant that does not appear to have 
been noted previously. Take $G=\HH$, the hexagonal lattice with connective constant 
$\mu=\sqrt{2+\sqrt{2}}\approx 1.84776$, see \cite{ds}. The ensuing lattice $\wt{\HH}$ is illustrated
in Figure \ref{fig:newlatt}, and its connective constant $\wt{\mu}$ satisfies $\mu^{-2}=h(\wt{\mu}^{-1})$,
which may be solved to obtain $\wt{\mu}\approx 1.75056$.
\end{example}

The proofs of Theorems \ref{thm:main2}--\ref{sf}
are based on the generating function $Z_v(x)$ of SAWs defined
by
\begin{equation}\label{eq:gf}
Z_v(x)=\sum_{w\in\Sigma(v)}x^{|w|},
\end{equation}
where $\Sigma(v)$ is the set of all finite SAWs starting from a given vertex $v$, and $|w|$ is the length
of $w$. Viewed as a power series, $Z_v(x)$ has radius of convergence 
$1/\mu$.

The connective constant is altered 
by application of the Fisher transformation as described in Theorem \ref{thm:main2}(a).
Critical exponents, on the other hand, have values that are not altered.
The reader is referred to \cite{bdgs, ms} and the references therein
for general accounts of critical exponents for SAWs.
The three exponents that have received most attention in the study of SAWs are 
as follows.

We consider only the case of SAWs in finite-dimensional spaces, thus excluding,
for example, the hyperbolic space of \cite{MadW}.
Suppose for concreteness that there 
exists a periodic, locally finite embedding of $G$ into $\RR^d$ with $d\geq 2$, and no
such embedding into $\RR^{d-1}$. The case of 
general $G$ has not been studied extensively, and most attention has been paid to the
hypercubic lattice $\ZZ^d$. 

\subsection*{The exponent $\g$.}
It is believed (when $d \ne 4$) that the generic behaviour of
$\s_n(v)$ is given by:
\begin{equation}\label{ex1}
\sigma_n(v)\sim A_vn^{\g-1}\mu^n,\qquad \text{as}\ n\to\infty,\  v\in V,
\end{equation}
for constants $A_v>0$ and  $\g\in\RR$.
The value of the `critical exponent' $\g$ is believed to depend on $d$ only, 
and not further on the choice of
graph $G$. Furthermore, it is believed (and largely proved, see the account in
\cite{ms}) that
$\g = 1$ when $d \ge 4$.
In the borderline case $d=4$, \eqref{ex1} should hold with $\g=1$ and subject to  
the correction factor $(\log n)^{1/4}$. 

\subsection*{The exponent $\eta$.}
Let $v,w\in V$, and
\begin{equation*}
Z_{v,w}(x)=\sum_{n=0}^{\infty}\sigma_n(v,w)x^k,\qquad x>0,
\end{equation*}
where $\sigma_n(v,w)$ is the number of $n$-step SAWs with endpoints $v,w$. It is known under certain 
circumstances that the generating functions $Z_{v,w}$ have radius of convergence $\mu^{-1}$ (see \cite[Cor.\ 3.2.6]{ms}),
and it is believed that there exists an exponent $\eta$ and constants $A_v'>0$ such that
\begin{equation*}
Z_{v,w}(\mu^{-1})\sim A_v'd_G(v,w)^{-(d-2+\eta)},\qquad \text{as } d_G(v,w)\to\infty,
\end{equation*}
where $d_G(v,w)$ is the graph-distance between $v$ and $w$.
Furthermore, $\eta$ satisfies $\eta=0$ when $d \ge 4$.  

\subsection*{The exponent $\nu$.}
Let $\Sigma_n(v)$ be the set of $n$-step SAWs from $v$, and write $\langle\cdot\rangle_n^v$ for
expectation with respect to the uniform measure on $\Sigma_n(v)$. Let $\|\pi\|$ be the 
graph-distance between the endpoints of a SAW $\pi$. 
It is believed (when $d\neq 4$) that there exists an
exponent $\nu$ and constants $A_v''>0$, such that
\begin{equation*}
\langle\|\pi\|^2\rangle_n^v\sim A_v''n^{2\nu}, \qquad v\in V.
\end{equation*}
As above, this should hold for $d=4$ subject to the inclusion of the 
correction factor $(\log n)^{1/4}$. It is believed that
$\nu=\frac12$ when $d \ge 4$.

\smallskip

The three exponents $\g$, $\eta$, $\nu$
are believed to be related through the so-called \emph{Fisher relation}
\begin{equation*}
\gamma=\nu(2-\eta).
\end{equation*}

In \cite[Sect.\ 3]{GrL2}, reasonable definitions of the three exponents are presented,
none of which depend on the existence of embeddings into $\RR^d$. 
Furthermore, it is proved that the values of the
exponents are unchanged under the Fisher transformation.  

\section{Bounds for connective constants}
\label{sec:lb}

Let $G$ be an infinite, connected, $\De$-regular graph.
How large or small can $\mu(G)$ be?  It is trivial that
$\s_n(v) \le \De(\De-1)^{n-1}$, whence $\mu(G) \le \De -1$.
It is not difficult to prove the strict inequality
$\mu(G) < \De-1$ when $G$ is quasi-transitive and contains a cycle
(see \cite[Thm 4.2]{GrL1}).
Lower bounds are harder to obtain.

\begin{theorem} \cite[Thm 4.1]{GrL1} \label{thm:lower}
Let $\Delta\geq 2$, and let $G$ be an infinite, 
connected, $\De$-regular, vertex-transitive
graph. Then $\mu(G)\geq\sqrt{\De-1}$ if either
\begin{letlist}
\item $G$ is simple, or
\item $G$ is non-simple and $\De \le 4$.
\end{letlist}
\end{theorem}

Note that, for the bridge graph $\BB_\De$ with $\De \ge 2$, we have the 
equality $\mu(\BB_\De)= \sqrt{\De-1}$. Theorem \ref{thm:lower}(a) answers
a question of Itai Benjamini.

\begin{question}\label{qn:q3}
What is the best universal lower bound in case {\rm(a)} above? In particular,
could it be the case that $\mu(G) \ge \phi$ for any infinite,
connected, cubic, vertex-transitive, \emph{simple} graph $G$?
\end{question}

\begin{question}\label{qn:q4}
Is it the case that $\mu(G) \ge \sqrt{\De-1}$
in the non-simple case {\rm(b)} with $\De > 4$?
\end{question}

Here is an outline of the proof of Theorem \ref{thm:lower}.
A SAW is called \emph{forward-extendable} if it is the initial
segment of some infinite SAW. 
Let $\siF_n(v)$ be the number of forward-extendable SAWs
starting at $v$.
Theorem \ref{thm:lower} is proved by showing as follows that 
\begin{equation}\label{g1}
\siF_{2n}(v) \ge (\De-1)^n.
\end{equation}
Let $\pi$ be a (finite) SAW from $v$, with final endpoint $w$. For a vertex $x \in \pi$ 
satisfying $x \ne w$, and an edge 
$e\notin\pi$ incident to $x$, the pair $(x,e)$
is called $\pi$-\emph{extendable}
if there exists an infinite SAW starting at $v$ whose initial segment
traverses $\pi$ until $x$, and then traverses $e$. 

First, it is proved subject to a 
certain condition $\Pi$ that, for any $2n$-step forward-extendable SAW  $\pi$,
there are at least $n(\Delta-2)$ $\pi$-extendable pairs. 
Inequality \eqref{g1} may be deduced from this statement.

The second part of the proof is to show that graphs
satisfying either (a) or (b) of the theorem satisfy 
condition $\Pi$. It is fairly simple to
show that (b) suffices, and it may well be reasonable to
extend the conclusion to values of $\De$ greater than $4$.

The growth rate $\muF$ of the number of forward-extendable SAWs has been studied 
further by Grimmett, Holroyd, and Peres \cite{GHLP}.
They show that $\muF = \mu$ for any infinite, connected, quasi-transitive
graph, with further results involving the numbers of
backward-extendable and doubly-extendable SAWs. 

\section{Strict inequalities for connective constants}\label{sec:si}

\subsection{Outline of results}
Consider a probabilistic model on a graph $G$, such as a percolation or random-cluster model 
(see \cite{G-rcm}). There is a parameter (perhaps `density' $p$ or `temperature' $T$)
and a `critical point' (usually written $\pc$ or $\Tc$). The numerical value
of the critical point depends  on the choice of graph $G$.
It is often important to understand whether a systematic change in the graph
causes a \emph{strict} change in the value of the critical point.
A general approach to this issue was presented by Aizenman and Grimmett \cite{AG} 
and developed further in \cite{BGK,G94} and \cite[Chap.\ 3]{G99}.
The purpose of this section is to review work of \cite{GrL3} directed
at the corresponding question for self-avoiding walks.

Let $G$ be a subgraph of $G'$, and suppose each graph is quasi-transitive.
It is trivial that $\mu(G) \le \mu(G')$. Under what conditions does 
the strict inequality $\mu(G) < \mu(G')$ hold? Two sufficient
conditions for the strict inequality 
are presented in \cite{GrL3}, and are reviewed here.
This is followed in Section \ref{sec:cayley} with a summary 
of the consequences for Cayley graphs.

\subsection{Quotient graphs}\label{sec:quot}
Let $G=(V,E)$ be a vertex-transitive graph. 
Let $\Ga$ be a subgroup of the automorphism group $\Aut(G)$ that
acts transitively, and let $\sA$ be a normal subgroup of $\Ga$
(we shall discuss the non-normal case later). There are several ways
of constructing a quotient graph $G/\sA$, the strongest of which (for our purposes) is
given next.
The set of neighbours of a vertex $v\in V$ is denoted by $\pd v$.

We denote by $\vG =(\olV,\vE)$ the \emph{directed} quotient graph
$G/\sA$ constructed as follows. Let $\approx$ be the equivalence relation on $V$ given by
$v_1 \approx v_2$ if and only if there exists
$\a\in\sA$ with $\a v_1=v_2$. The vertex-set $\olV$
comprises the equivalence classes of $(V,\approx)$, that is, the orbits $\ol v := \sA v$
as $v$ ranges over $V$. For $v,w \in V$, we place $|\pd v \cap \ol w|$ directed edges
from $\ol v$ to $\ol w$ (if $\ol v = \ol w$,
these edges are directed loops).

\begin{example}
As a simple example of a quotient graph, consider the square lattice $G=\ZZ^2$ and let $m \ge 1$.
Let $\Ga$ be the set of translations of $\ZZ^2$, and let $\sA$ be the 
normal subgroup of $\Ga$ generated by the map that sends each $(i,j)$ to
$(i+m,j)$. The quotient graph $G/\sA$ is the square lattice `wrapped
around a cylinder', with each edge replaced by two oppositely directed edges.

A second example is presented using the language of Cayley graphs in Example \ref{solattice}.
\end{example}

Since $\vG$ is obtained from $G$
by a process of identification of vertices and edges, it is natural to ask whether 
$\mu(\vG) < \mu(G)$. Sufficient conditions for this strict inequality are presented next. 

Let $L=L(G,\sA)$ be the length of the shortest SAW of $G$
with (distinct) endpoints in the same orbit.
Thus, for example, $L=1$ if $\vG$ possesses a directed loop.
A group is called \emph{trivial} if it comprises the identity only.

\begin{theorem} \cite[Thm 3.8]{GrL3} \label{si}
Let $\Ga$ act transitively on $G$, and let
$\sA$ be a non-trivial, normal subgroup of $\Ga$. 
The connective constant $\vmu=\mu(\vG)$ satisfies $\vmu < \mu(G)$ if: either
\begin{letlist}
\item $L \ne 2$, or
\item $L=2$ and either of the following holds:
\begin{romlist}
\item $G$ contains some $2$-step
 SAW $v\ (=w_0),w_1,w_2\ (=v')$ satisfying $\ol v = \ol v'$ and $|\pd v \cap \ol w_1| \ge 2$, 
\item $G$ contains some SAW $v\,(=w_0),w_1,w_2,\dots,w_l\,(=v')$ satisfying 
$\ol v = \ol v'$, $\ol w_i \ne \ol w_j$ for $0 \le i < j < l$, and furthermore $v' = \a v$ for
some $\a \in \sA$ which fixes no $w_i$.
\end{romlist}
\end{letlist}
\end{theorem} 

\begin{question}\label{qn:strict}
In the situation of Theorem \ref{si}, can one \emph{calculate} an explicit $\eps=\eps(G,\sA)>0$
such that $\mu(G) - \mu(\vG) > \eps$? A partial answer is provided at \cite[Thm 3.11]{GrL3}.
\end{question}

We call $\sA$ \emph{symmetric} if 
\begin{equation*}
|\partial v\cap\ol{w}|=|\partial w\cap\ol{v}|,\qquad v,w\in V.
\end{equation*}
Consider the special case $L=2$ of Theorem \ref{si}.
Condition (i) of Theorem \ref{si}(b)
holds if $\sA$ is symmetric, since $|\pd w \cap \ol v| \ge 2$.
Symmetry of $\sA$ is implied by unimodularity, but we prefer to
avoid the topic of unimodularity in this review, referring the 
reader instead to \cite[Sect.\ 3.5]{GrL3} or \cite[Sect.\ 8.2]{LyP}.

\begin{example}
Conditions (i)--(ii) of Theorem \ref{si}(b) are necessary in the case $L=2$,
in the sense illustrated by the following example. Let $G$ be the 
infinite cubic tree with a distinguished end $\om$.
Let $\Ga$ be the set of automorphisms
that preserve $\omega$, and let $\sA$ be the normal subgroup
generated by the interchange of two children of a given vertex $v$
(and the associated relabelling of their descendants). 
The graph $\vG$ is isomorphic to that obtained from $\ZZ$ by
replacing each edge by two directed edges in one direction and one in the reverse direction.
It is easily seen that $L=2$, but that
neither (i) nor (ii) holds.
Indeed, $\mu(\vG) = \mu(G) = 2$.
\end{example}

\begin{figure}[htbp]
\centering
\includegraphics*{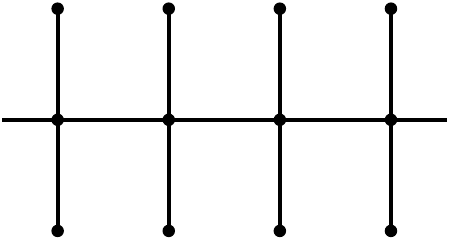}
\caption{The pattern is extended infinitely in both directions.}
\label{fig:reflect}
\end{figure}

The conclusion of Theorem \ref{si} is
generally invalid under the weaker assumption that $\sA$ acts quasi-transitively on $G$.  
Consider, for example, the graph $G$ of Figure \ref{fig:reflect},
with $\sA = \{\id,\rho\}$ where $\rho$
is reflection in the horizontal axis.
Both $G$ and its quotient graph have connective constant $1$.

The proof of Theorem \ref{si} follows partly the general approach of Kesten
in his pattern theorem, see \cite{hkI} and \cite[Sect.\ 7.2]{ms}.
Any $n$-step
SAW $\vpi$ in the directed graph $\vG$  lifts to a SAW $\pi$ in the larger
graph $G$. The
idea is to show that `most' such $\vpi$ contain at least $an$ sub-SAWs 
for which the corresponding sub-walks of $\pi$ may be replaced by SAWs
on $G$.
Different subsets of these sub-SAWs of $\vG$ give rise to different SAWs 
on $G$. The number of such subsets grows exponentially in $n$, and this
introduces an exponential `entropic' factor in the counts of SAWs.

Unlike Kesten's proof and its subsequent 
elaborations by others, our results
apply in the general setting of vertex-transitive graphs, and they
utilize algebraic and combinatorial techniques.

We discuss next the assumption of normality of $\sA$ in
Theorem \ref{si}. The (undirected) simple quotient graph $\olG=(\olV,\olE)$ may be defined as follows
even if $\sA$ is not a normal subgroup of $\Ga$. As before, the vertex-set $\olV$ is the set
of orbits of $V$ under $\sA$. Two distinct orbits $\sA v$, $\sA w$ are declared adjacent 
in $\olG$ if there exist $v' \in \sA v$ and $w' \in \sA w$ with $\langle v',w'\rangle \in E$.
We write $\olG = G_\sA$ to emphasize the role of $\sA$.

The relationship between the \emph{site percolation} critical points of $G$ and 
$G_\sA$ is the topic of
a conjecture of Benjamini and Schramm \cite{BenS96}, which
appears to make the additional assumption that $\sA$
acts freely on $V$.  The last assumption is stronger
than the assumption of unimodularity.

We ask for an example in which the non-normal case is essentially different from the
normal case.

\begin{question}\label{qn:q5}
Let $\Ga$ be a subgroup of $\Aut(G)$ acting transitively on $G$. 
Can there exist a non-normal subgroup $\sA$
of $\Ga$ such that: {\rm(i)} the quotient graph $G_\sA$ is vertex-transitive,
and {\rm(ii)} there exists no normal subgroup $\sN$ of some
transitively acting $\Ga'$ such that
$G_\sA$ is isomorphic to $G_\sN$? 
Might it be relevant to assume that $\sA$ acts freely on $V$?
\end{question}

We return to connective constants with the following question.

\begin{question}\label{qn:qn7'}
Is it the case that $\mu(G_\sA) < \mu(G)$ under the assumption that 
$\sA$ is a non-trivial (not necessarily normal) subgroup of $\Ga$ 
acting freely on $V$, such that $G_\sA$ is vertex-transitive?
\end{question}

The proof of Theorem \ref{si}  may be adapted
to give an affirmative answer to Question \ref{qn:qn7'} subject a certain extra condition on $\sA$,
see \cite[Thm 3.12]{GrL3}. Namely, it suffices that there exists $l \in \NN$
such that $G_\sA$ possesses a cycle of length $l$ but $G$ has no cycle of this length.

\subsection{Quasi-transitive augmentations} \label{sec:qta}
We consider next the systematic addition of new edges, and the effect thereof
on the connective constant.
Let $G=(V,E)$ be an infinite, connected, vertex-transitive, simple
graph. From $G$, we derive a second graph
$\ollG=(V,\ollE)$ 
by adding further edges to $E$, possibly in parallel to existing edges. 
\emph{We assume that $E$ is a proper subset of $\ollE$}, and introduce
next a certain technical property. Let $\Ga$ be a subgroup of $\Aut(G)$ that acts transitively.

\begin{definition}\label{def:fcp}
A subgroup $\sA$ of $\Ga$ is said to
\emph{have the \fcp} (relative to $\Ga$) with \emph{root} $\rho \in V$ if there exist
$\nu_0,\nu_1,\dots,\nu_s \in \Ga$, with $\nu_0=\id$ and $s<\oo$, such
that $V$ is partitioned as $\bigcup_{i=0}^s \nu_i \sA \rho$.
It is said simply to have the \fcp\ if it has this property with some root.
\end{definition}

This definition is somewhat technical. The principal situation in which 
$\sA$ has the \fcp\ arises,
as stated in Theorem \ref{thm:sufft}, when $\sA$ is a normal subgroup of $\Ga$
that acts quasi-transitively on $G$.

\begin{theorem} \cite[Thm 3.2]{GrL3}\label{qta}
Let $\Ga$ act transitively on $G$, and let $\sA$ be a subgroup of $\Ga$ 
with the finite coset property. If $\sA\subseteq\Aut(\ollG)$,
then $\mu(G)<\mu(\ollG)$.
\end{theorem}

\begin{example} \label{ex:sqt}
Let $\ZZ^2$ be the square lattice, with $\sA$  the group
of its translations.  The triangular lattice $\TT$ is obtained from $\ZZ^2$ by adding the edge 
$e=\langle 0,(1,1)\rangle$
together with its images under $\sA$, where $0$ denotes the origin. 
Since $\sA$ is a normal subgroup
of itself with the \fcp, it follows that $\mu(\ZZ^2) < \mu(\TT)$. This example
may be extended to augmentions by other periodic families of new edges,
as explained in \cite[Example 3.4]{GrL3}.
\end{example}

\begin{question}\label{qn:strict2}
In the situation of Theorem \ref{qta}, 
can one calculate $\eps>0$ such that $\mu(\ollG)-\mu(G) > \eps$?
(See the related  Question \ref{qn:strict}.)
\end{question}

Two classes of subgroup $\sA$ with the finite coset property are given as follows.

\begin{theorem} \cite[Prop.\ 3.3]{GrL3}\label{thm:sufft}
Let $\Ga$ act transitively on $G$, and $\rho\in V$.
The subgroup $\sA$ of $\Ga$ has the finite coset property
with root $\rho$ if either of the following holds.
\begin{numlist}
\item $\sA$ is a normal subgroup of $\Ga$ which acts quasi-transitively on $G$.
\item The index $[\Ga:A]$ is finite.
\end{numlist}
\end{theorem}

\begin{figure}[htbp]
\centerline{\includegraphics[width=0.7\textwidth]{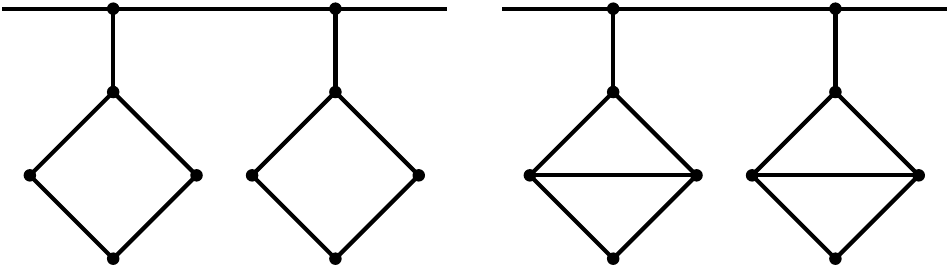}}
   \caption{The pair $G$, $\ollG$. The graphs are extended in both directions. 
Each graph is quasi-transitive with connective constant $1$, 
and the second is obtained
from the first by a quasi-transitive addition of edges.}
\label{fig:qusi}
\end{figure}

It would be insufficient to assume only quasi-transitivity in Theorem \ref{qta}. 
Consider, for example, the pair
$G$, $\ollG$ of Figure \ref{fig:qusi}, each of which has connective constant 1.

By Theorem \ref{thm:sufft}(a), $\mu(G) < \mu(\ollG)$ if
$\sA$ is a normal subgroup of some transitive $\Ga$, and $\sA$ acts quasi-transitively.
Can we dispense with the assumption of normality?

\begin{question}\label{qn:nonnormal}
Let $\Ga$ act transitively on $G$, and let $\sA$ be a subgroup of $\Ga$
that acts quasi-transitively on $G$. 
If $\sA \subseteq \Aut(\ollG)$, is it necessarily the case that
$\mu(G) < \mu(\ollG)$?
\end{question}

A positive answer would be implied by an affirmative answer to the following question.
For $\sA \subseteq \Ga \subseteq \Aut(G)$, we say that \emph{$\Ga\setminus \sA$ acts freely on
$V$} if: whenever $\g \in \Ga$ and $v \in V$ satisfy $\g v = v$, then $\g \in \sA$. 

\begin{question}\label{qn:nonnormal2}
Let $G$ be a vertex-transitive graph, and let $\sA$ be a subgroup of $\Aut(G)$
that acts quasi-transitively on $G$. Does there exist (or, weaker, \emph{when} does there exist)
a subgroup $\Ga$ of $\Aut(G)$ acting transitively on $G$ such that $\sA \subseteq \Ga$ and
$\Ga \setminus \sA$ acts freely on $V$?
\end{question}

See \cite[Prop.\ 3.6]{GrL3} and the further discussion therein.

\section{Connective constants of Cayley graphs}\label{sec:cayley}

Theorems \ref{si}--\ref{qta} have the following implications for Cayley graphs.
Let $\Ga$ be an infinite group with a finite generator-set $S$, 
where $S$ satisfies $S=S^{-1}$ and $\id \notin S$. Thus, $\Ga$ has a presentation as
$\Ga = \langle S \mid R\rangle$ where $R$ is a set of relators.
The Cayley graph $G=G(\Ga,S)$ is defined as follows.
The vertex-set $V$ of $G$ is the set of elements of $\Ga$. 
Distinct elements $g,h \in V$
are connected by an edge if and only if there exists $s \in S$ such that $h=gs$.
It is easily seen that $G$ is connected and vertex-transitive, and it is standard that
$G$ is unimodular and hence symmetric. Accounts of
Cayley graphs may be found in \cite{bab95}  and \cite[Sect.\ 3.4]{LyP}. 

Let $s_1s_2\cdots s_l = \id$ be a relation. 
This relation corresponds to 
the closed walk $(\id,s_1, s_1s_2,\dots, s_1 s_2\cdots s_l=\id)$
of $G$ passing through the identity $\id$. 
Consider now the effect of adding a further relator.
Let $t_1,t_2,\dots,t_l \in S$ be such that 
$\rho := t_1t_2\cdots t_l$ satisfies $\rho \ne \id$, and write
$\Ga_\rho = \langle S \mid R\cup\{\rho\}\rangle$.
Then $\Ga_\rho$ is isomorphic to the quotient group $\Ga/\sN$ where
$\sN$ is the normal subgroup of $\Ga$ generated by $\rho$.

\begin{theorem} \cite[Corollaries 4.1, 4.3]{GrL3}\label{caleydecrease}
Let $G=G(\Ga,S)$ be the Cayley graph of
the infinite, finitely generated group $\Ga = \langle S\mid R\rangle$.

\begin{letlist}
\item
Let $G_\rho= G(\Ga_\rho, S)$ be the Cayley graph obtained by adding to $R$ a further non-trivial relator $\rho$. 
Then $\mu(G_\rho) <\mu(G)$.

\item
Let $w \in\Ga$ satisfy $w \ne \id$, $w \notin S$,
and let $\olG_w$ be the Cayley graph of the group obtained by adding $w$
(and $w^{-1}$) to $S$. Then $\mu(G) < \mu(\olG_w)$.
\end{letlist}
\end{theorem} 

\begin{example}\label{solattice}
The \emph{square/octagon} lattice, otherwise known as the Archimedean
lattice $(4,8^2)$, is the 
Cayley graph of the group with generator set 
$S=\{s_1,s_2,s_3\}$ and relators 
$$
\{s_1^2, s_2^2, s_3^2, 
s_1s_2s_1s_2,s_1s_3s_2s_3s_1s_3s_2s_3\}.
$$
(See \cite[Fig.\ 3]{GrL3}.)
Adding the further
relator $s_2s_3s_2s_3$, we obtain a graph isomorphic
to the ladder graph of
Figure \ref{fig:ladder-hex}, whose connective constant is the golden mean $\phi:=\frac12(\sqrt{5}+1)$. 

By Corollary \ref{caleydecrease}(a),
the connective constant $\mu$ of the square/octagon lattice is strictly greater than 
$\phi = 1.618\dots$. The best lower bound currently known appears to
be $\mu > 1.804\dots$, see \cite{j04}.
\end{example}

\begin{example}
The square lattice $\ZZ^2$ is the Cayley graph
of the abelian group with $S = \{a,b\}$
and $R=\{aba^{-1}b^{-1}\}$. We add a generator $ab$
(and its inverse), thus adding a diagonal to each square of $\ZZ^2$.
Theorem \ref{caleydecrease}(b) implies the standard inequality $\mu(\ZZ^2)< \mu(\TT)$
of Example \ref{ex:sqt}.
\end{example}

\section*{Acknowledgements}
This work was supported in part by the Engineering
and Physical Sciences Research Council under grant EP/103372X/1.
GRG acknowledges valuable conversations
with Alexander Holroyd concerning Questions \ref{qn:q5} and 
\ref{qn:nonnormal2}. The authors are grateful to the referee for 
a careful reading and numerous suggestions.

\bibliography{saw-final2}
\bibliographystyle{amsplain}

\end{document}